\documentclass[doublespacing]{elsart}

 \usepackage{graphics}
 \usepackage{graphicx}
 \usepackage{epsfig}

\usepackage{amssymb}
\usepackage{amsmath}
\usepackage{subfigure}
\begin{document}

\begin{frontmatter}



\title{Actuator and sensor placement in linear advection PDE with building systems application}


\author[vaidya]{U. Vaidya\corauthref{vaidya1}}
\thanks[vaidya1]{Corr. author Email: ugvaidya@iastate.edu}
\author[rajeev]{R. Rajaram}
\author[dasgupta]{and S. Dasgupta}
\address[rajeev]{\scriptsize Department of Math. Sci., 3300, Lake Rd West, Kent State University, Ashtabula, OH - 44004}
\address[vaidya]{\scriptsize Dept. of Elec. and Comp. Engg.,
Iowa State University, Ames, IA 50011, Email: ugvaidya@iastate.edu}
\address[dasgupta]{\scriptsize Dept. of Elec. and Comp. Engg.,
Iowa State University, Ames, IA 50011, Email: dasgupta.sambarta@gmail.com}

\begin{abstract}
We study the problem of actuator and sensor placement in a linear advection partial differential equation (PDE). The problem is motivated by its application to actuator and sensor placement in building systems for the control and detection of a scalar quantity such as temperature and contaminants.
We propose a gramian based approach to the problem of actuator and sensor placement. The special structure of the advection PDE is exploited to provide an explicit formula for the controllability and observability gramian in the form of a multiplication operator. The explicit formula for the gramian, as a function of actuator and sensor location, is used to provide test criteria for the suitability of a given sensor and actuator location. Furthermore, the solution obtained using gramian based criteria is interpreted in terms of the flow of the advective vector field. In particular, the almost everywhere stability property of the advective vector field is shown to play a crucial role in deciding the location of actuators and sensors. Simulation results are performed to support the main results of this paper.
\end{abstract}
\end{frontmatter}

% main text
\section{Introduction}
In this paper, we study the problem of actuator and sensor placement in a linear advection partial differential equation. The problem is motivated by its application to actuator and sensor location in building systems for the purpose of control of temperature and detection of contaminants.
Building systems in US account for $39$ percent of total energy consumption \cite{building_system_website}. Design of efficient building systems not only has a significant economic benefit but also  social and environmental benefits. Social benefits arise due to improved overall quality of life by enhancing occupant health, comfort and heightened aesthetic qualities. Improvement in water and air quality and reduced waste lead to environmental benefits.

The optimal placement of actuators and sensors in a building system is a difficult problem due to the complex physics that is involved. The governing equations for building system fluid flows and scalar densities are coupled nonlinear partial differential equations subjected to disturbances, various sources of uncertainties, and complicated geometry. Analysis of the building system  with its full scale complexity leads to a finite element based computational approach to the actuator and sensor placement problem \cite{optimal_actuator_sensor_cfd}. Such a purely computational based approach provides little insight into the obtained solution. An alternate system theoretic and dynamical systems based approach under some simplifying assumptions and physics can also be pursued \cite{surana_building,surana_building_dynamical_systems}. Such an approach provides useful insights and guidelines to the complex control problems involved in building system applications. In this paper, we pursue a similar approach for the location of actuators and sensors problem in a building system.

Under some simplifying assumptions and physics \cite{surana_building, surana_building_dynamical_systems}, the system equations are modeled in the form of a linear advection partial differential equation with inputs and outputs. We propose a gramian based approach to the actuator and sensor location problem.
%Because of the simplifying assumptions made in the modeling of system equations, the results in the paper cannot be directly applied to the building system application.
The results are an important first step towards its application to building systems. However, further research needs to be done for relaxing some of the simplifying assumptions made in this paper for the applicability of these results for the building systems problem. We believe that the analytical methods developed in this paper combined with computational techniques involving detailed physics of building systems is a right approach moving forward.  The main contribution of this paper is in providing explicit formula for the controllability and observability gramians as a function of actuator and sensor locations and the advection velocity field. These explicit formulas for the gramians are used to provide test criteria for deciding the location of sensors and actuators. Technical conditions for the existence of infinite time gramians are also provided. In particular, we prove that the infinite time controllability and observability gramians are well defined for almost everywhere stable and asymptotically stable advection vector fields respectively. We provide simulation results using a two dimensional fluid flow vector field for the computation of the finite time controllability and observability gramian.

An excellent review and classification of sensor and controller positioning for distributed parameter systems can be found in \cite{kubrusly}, where most of the methods involve a finite dimensional approximation of the infinite dimensional system, either before or after solving an optimization problem using the point spectrum
of the infinitesimal generator. It was noted in \cite{kubrusly} that such an approximation based method will not work for wave type systems because of the finite speed of propagation. \cite{eljai} is an excellent book on sensor and actuator placement for distributed parameter systems governed by heat and diffusive type processes. \cite{gawronski} describes sensor and  actuator placement for flexible structures. A combinatorial optimization approach for linear time invariant systems
based on integer programming using the controllability and observability gramians for sensor and actuator placement can be found in \cite{georges}.

The linear advection equation considered in this paper is akin to a unidirectional wave equation whose wave speed is governed by a nonlinear smooth vector field $f(x)$. Hence, actuator and sensor placement analysis based on a finite dimensional approximation as developed in earlier references will not work well for our problem. Our selection criteria for actuators
and sensors uses the idea of controllability and observability gramians, but differs from what is seen in the literature slightly. The advection equation has a fundamental limitation for control as described in Theorem (\ref{pde_controllability}), in the sense that placing an actuator on the set $B$ can
only affect states $\rho$ whose support is $\mathcal{R}_B^{\tau} = \cup_{t=0}^{\tau} \phi_t(B)$. This set $\mathcal{R}_B^{\tau}$
is precisely the support of the controllability gramian $\mathcal{C}_B^{\tau}$ for the advection equation (see Claim (\ref{claimcontrollabilitygramian})).
This is the main reason why we consider choosing a set $B$ for actuators that maximizes the support of $\mathcal{C}_B^{\tau}$. For sets $B$ that give
the same support for $\mathcal{C}_B^{\tau}$, we choose the one that gives lesser $L^2$ norm, since that will minimize the control effort (see Theorem
(\ref{pde_controllability}), where the minimum norm control formula has the controllability gramian appearing in the denominator).

The organization of the paper is as follows. In section \ref{preliminary}, we describe the problem and some preliminaries from the theory of partial differential equations. In section \ref{main}, we present the main results of the paper. In section \ref{Infinitetimegramians}, we discuss technical conditions for the existence of the infinite time controllability and observability
gramians. Simulation results are presented in section \ref{simulation} followed by conclusion in section \ref{conclusion}.
\section{Preliminaries}\label{preliminary}
We study the problem of optimal location of actuator in a linear advection partial differential equation. The motivation for this problem comes from the optimal location of actuators for the control of a scalar quantity, such as temperature or contaminants, in a room denoted by $\rho(x,t)$.

In building system applications, the evolution of $\rho(x,t)$, is governed by the velocity field $v(x,t)$ of the fluid flow. This velocity field is obtained as a solution of the following Navier Stokes equation:
\begin{eqnarray}
&\frac{\partial v(x,t)}{\partial t}+v(x,t)\cdot \nabla v(x,t)=-\nabla p(x,t)+\frac{1}{Re}\bigtriangleup v(x,t)\nonumber\\
& \nabla \cdot v(x,t)=0,\label{NS}
\end{eqnarray}
where $x\in X \subset \mathbb{R}^N$ (with $N=2$ or $3$) is the domain of the room, $v(x,t)$ is the velocity field, $p(x,t)$ is the pressure, and $Re$ is the Reynolds number. The evolution of the scalar quantity $\rho(x,t)$ is governed by the following linear controlled partial differential equation
\begin{eqnarray}
\frac{\partial \rho}{\partial t}+v(x,t)\cdot \nabla \rho(x,t) &=&\frac{1}{Pr Re} \bigtriangleup \rho(x,t)+\sum_{k=1}^N\chi_{B_k}(x)u_k(t)\nonumber\\
y_k(x,t)&=&\chi_{A_k}(x)\rho(x,t),\;\;\;k=1,\ldots,M\label{advection_diffusion}
\end{eqnarray}
where $Pr$ is the Prandtl number, $\chi_{A_k}(x)$ is the indicator function on set $A_k\subset X$, and $u_k(t)\in \mathbb{R}$ is the control input for $k=1,\ldots,N$. The form of control input $\chi_B(x)u(t)$ and output measurement $\chi_A(x)\rho(x,t)$ is motivated by the fact that the actuation and sensing can be exercised only on a small region $B$ and $A$ of the physical space $X$ respectively.
\begin{rem}
\label{outputequation}
The form of output equation in (\ref{advection_diffusion}) is different than the one usually considered in the literature,
where the sensors have access to the average state information on a set (i.e., $y_k(t)=\int_{A_k} c_k(x)\rho(x,t)$). The interpretation in our case is that
the sensors have pointwise state information from the sets $A_k$. We choose
the form in (\ref{advection_diffusion}) because it allows us to compute the observability gramians as a explicit function of section location set $A_k$. Our proposed approach can also be applied to the case where the sensors have access to average state information, however, the observability gramian in that case will be a complicated function of the sensor location set $A_k$. Furthermore this form of output measurement is also dual to the input actuation term, in particular to Eq. (\ref{transportcontrol}).
\end{rem}

The objective is to determine the optimal location of actuators and sensors, and hence the determination of indicator function $\chi_{B_k}(x)$ and $\chi_{A_k}(x)$. The terms $v(x,t)\cdot \nabla \rho(x,t)$ and $ \bigtriangleup \rho(x,t)$ in (\ref{advection_diffusion}) correspond to advection and diffusion respectively, with $D=\frac{1}{RePr}$ being the diffusion constant. Note that the advection diffusion equation (\ref{advection_diffusion})  is decoupled from the Navier Stokes equation (\ref{NS}). In the case where the scalar density is temperature, this decoupling corresponds to the assumption that buoyancy forces have negligible or no effect. Furthermore, for simplicity of presentation of the main results of this paper, we now make the following assumptions.
\begin{assum}\label{assume_steady} We replace the time varying velocity field  $v(x,t)$ responsible for the advection of scalar density with the mean velocity field $f(x)$ i.e.,
\[f(x):=\frac{1}{T}\int_0^T v(x,t)dt\]
\end{assum}
\begin{rem}
Typically the velocity field information $v(x,t)$ is available over a finite time interval $[0,T]$ either from a simulation or from an experiment. %In such a case, the steady velocity field can be obtained using the time average of $v(x,t)$ i.e.,
%\[f(x):=\frac{1}{T}\int_0^T v(x,t)dt\]
Assumption (\ref{assume_steady}) corresponds to linearizing the linear advection PDE along the mean flow field $f(x)$.
It follows that if $v(x,t)$ is volume preserving i.e., $\nabla \cdot v(x,t)=0$, then $\nabla \cdot f(x)=0$ as well.
\end{rem}
\begin{assum}\label{assume_zero_diffusion}
Again for simplicity of presentation of the main results of this paper, we assume that the diffusion constant $D$ in the advection diffusion equation (\ref{advection_diffusion}) is zero. As we see in the simulation section, the assumption of zero diffusion constant is justified.
\end{assum}
%\begin{rem}
%As we see in the simulation section, the assumption of zero diffusion constant is justified.
%\end{rem}
We next discuss a few preliminaries on semigroup theory of partial differential equations.
Consider the following ordinary differential equation (ODE):
\begin{eqnarray}
\dot x=f(x),\;\;\;\;\;x(0)=x_0,\label{ode}
\end{eqnarray}
where $x\in X\subset \mathbb{R}^N$ a compact set. We denote by $\phi_t(x)$ the solution of ODE (\ref{ode}) starting from the initial condition $x$.
ODE (\ref{ode}) is used to define two linear infinitesimal operators, ${\cal A}_K:L^2(X)\to L^2(X)$ and ${\cal A}_{PF}: L^2(X)\to L^2(X)$ defined as follows:
\[{\cal A}_K \rho =f\cdot \nabla \rho,\;\;\;\; {\cal A}_{PF}\rho=-\nabla \cdot (f\rho).\]
The domains of the above operators are given as follows:
\[D({\cal A}_K ) = \{\rho \in H^1(X): \rho|_{\Gamma_o} = 0\},\]\[D({\cal A}_{PF}) = \{\rho \in H^1(X): \rho|_{\Gamma_i} = 0\},\]
where $\Gamma_o$ and $\Gamma_i$ are the outflow and inflow portions of the boundary $\partial X$ defined as follows:
\[\Gamma_o = \{x \in \partial X: f \cdot \eta > 0\},\;\;\;
\Gamma_i = \{x \in \partial X: f \cdot \eta < 0\},\]
where $\eta$ is the outward normal to the boundary $\partial X$.
The semigroups corresponding to the ${\cal A}_K$ and ${\cal A}_{PF}$ are called as Koopman $(\mathbb{U}_t)$ and Perron-Frobenius $(\mathbb{P}_t)$ operators respectively. These operators are defined as follows:
\[\mathbb{U}_t:L^2(X)\to L^2(X),\;\; (\mathbb{U}_t\rho)(x)=\rho(\phi_t(x)),\]
\[\mathbb{P}_t:L^2(X)\to L^2(X),\;\; (\mathbb{P}_t\rho)(x)=\rho(\phi_{-t}(x))\left|\frac{\partial \phi_t(x)}{\partial x}\right|^{-1},\]
where $|\cdot|$ denotes the determinant. These semigroups can be shown to satisfy the following partial differential equations \cite{Lasota}:
\[\frac{\partial \rho}{\partial t}-{\cal A}_K\rho=0,\rho|_{\Gamma_o}=0;\;\;\;\;\;
\frac{\partial \rho}{\partial t}-{\cal A}_{PF}\rho=0,\rho|_{\Gamma_i}=0.\]
 The Koopman and Perron-Frobenius semigroup operators and their infinitesimal generators are adjoint to each other i.e.,
\[\int_X (\mathbb{P}_t \rho_1)(x)\rho_2(x)d x=\int_X \rho_1(x)(\mathbb{U}_t \rho_2)(x)dx ~\forall \rho_1,\rho_2 \in L^2(X).\]

\section{Main results}\label{main}
The gramian based approach is one of the systematic approaches available for the optimal placement of actuators and sensors. Controllability and observability gramians measure the relative degree of controllability and observability of various states in the state space. Using the gramian based approach, actuators and sensors are placed at a location where the degree of controllability and observability of the least controllable and observable state is maximized \cite{Curtain,paganini_dull}.

\subsection{Controllability gramian}
For the construction of the controllability gramian, the advection-diffusion partial differential equation (\ref{advection_diffusion}) using assumptions (\ref{assume_steady}) and (\ref{assume_zero_diffusion}) for a single input case can be written as follows:
\begin{eqnarray}
\label{transportcontrol}
&&\frac{\partial \rho}{\partial t} + \nabla \cdot (f(x)\rho) = \chi_B(x) u(x,t);\;\;\;
\\\nonumber&&\rho|_{\Gamma_i} =  0;\;\; \rho(x,0) = \rho_0(x).
\end{eqnarray}
In Eq. (\ref{transportcontrol}), we are assuming that the control input $u$ is both a function of spatial variable $x$ and time $t$. This assumption will typically not be satisfied in the building system application, however, making this assumption allows us to use existing results from linear PDE theory in the development of controllability gramian \cite{Curtain}. Furthermore, since $m(X)>>m(B)$, where $m$ is the Lebesgue measure, we expect the main conclusions of this paper to hold even when $u$ is assumed to be only a function of time.
The set $B$ is the region of control in the state space $X$, and $u(x,t) \in L^2([0,\tau]:L^2(B))$
i.e., we have a control input that is square integrable in time and space, acting on the set $B$.
The solution to (\ref{transportcontrol}) is given by the following:

\begin{eqnarray*}
\rho(x,t) &&= \mathbb{P}_t\rho_0(x) + \int_0^t \mathbb{P}_{t-s} (\chi_B(x) u(x,s)) ds.
\end{eqnarray*}
We define the controllability operator $\mathbf{\mathcal{B}}^{\tau}:L^2([0,\tau]:L^2(B)) \rightarrow L^2(X)$ as follows:
\begin{eqnarray}
\label{controllabilityoperator}
(\mathbf{\mathcal{B}}^{\tau} u)(x) &&:= \int_0^t \mathbb{P}_{t-s} (\chi_B(x) u(x,s)) ds.
\end{eqnarray}
The adjoint of the controllability operator $\mathbf{\mathcal{B}}^{\tau *}:L^2(X) \rightarrow L^2([0,\tau]:L^2(B))$
can be calculated and is given as follows:
\begin{eqnarray}
\label{adjointofcontrollabilityoperator}
(\mathbf{\mathcal{B}}^{\tau *} z)(x,s) = \chi_B(x) \mathbb{U}_{(\tau -s)} z(x).
\end{eqnarray}
We have the following theorem on the controllability property of the PDE (\ref{transportcontrol}).
\begin{thm}\label{pde_controllability}
Let $\mathcal{R}^{\tau} = \cup_{t=0}^{\tau} \phi_t(B)$. The PDE (\ref{transportcontrol}) is exactly controllable in a given
time $\tau > 0$ for all initial and terminal states in the space $L^2(\mathcal{R}^{\tau})$ i.e. given initial and terminal states $\rho_0(x)$
and $\rho_{\tau}(x)$ in $\mathcal{R}^{\tau}$, there exists a control $u(x,t) \in L^2([0,\tau]:L^2(B))$ such that
$\rho(x,0)=\rho_0(x)$, and $\rho(x,\tau)= \rho_{\tau}(x)$, where $\rho(x,t)$ is the solution of (\ref{transportcontrol}).
\end{thm}
\begin{pf}
We prove the theorem by showing the following, which is equivalent to showing that the range of the controllability operator
$\mathbf{\mathcal{B}}^{\tau}$ is the same as $L^2(\mathcal{R}^{\tau})$:
\begin{enumerate}
\item \begin{equation}\label{approximatecontrollability} \mathbf{\mathcal{B}}^{\tau *} z = 0 ~\forall (x,s) \in B \times [0,\tau] \Rightarrow z = 0 \mbox{ in } L^2(\mathcal{R}^{\tau}).\end{equation}
\item The range of $\mathbf{\mathcal{B}}^{\tau}$ is closed.
\end{enumerate}
Assume that $\mathbf{\mathcal{B}}^{\tau *}z = \chi_B(x) \mathbb{U}_{(\tau -s)} z(x) = 0 ~\forall (x,s) \in B \times [0,\tau]$.
The assumption simply means that $z=0$ on $\cup_{t=-\tau}^0 \phi_t(B)$.
Since the set $\cup_{t=-\tau}^0 \phi_t(B)$ evolves into $\mathcal{R}^{\tau} = \cup_{t=0}^{\tau} \phi_t(B)$,
we have that $z=0$ in $L^2(\mathcal{R}^{\tau})$ .
Next, we recall that $\mathcal{B}^{\tau}:L^2([0,\tau]:L^2(B)) \rightarrow L^2(X)$
is defined by $(\mathcal{B}^{\tau} u)(x):= \int_0^{\tau} \mathbb{P}_{-(\tau - s)} \chi_B(x) u(x,s) ds$.
Let us assume that $u_n(x,t) \rightarrow u(x,t)$ is a convergent sequence in $L^2([0,\tau]:L^2(B))$.
We need to show that $(\mathcal{B}^{\tau} u_n)(x) \rightarrow (\mathcal{B}^{\tau} u)(x)$ in
$L^2(X)$. We have the following, where we have used $||\mathbb{P}_{t}||_{L^2(X)}\leq M_{\omega} e^{\omega t}$ from
the semigroup property of $\mathbb{P}_t$:
\begin{eqnarray}
\nonumber&& ||(\mathcal{B}^{\tau} u_n)(x) - (\mathcal{B}^{\tau} u)(x)||_{L^2(X)}^2
\\\nonumber&&= \int_X \int_0^{\tau} |\mathbb{P}_{(\tau - s)} \chi_B(x) (u_n(x,s) - u(x,s))|^2 ds dx
%\\\nonumber&& =  \int_0^{\tau} \int_X|\mathbb{P}_{(\tau - s)} \chi_B(x) (u_n(x,s) - u(x,s))|^2 dx ds
\\\nonumber&& = \int_0^{\tau} ||\mathbb{P}_{(\tau - s)} \chi_B(x) (u_n(x,s) - u(x,s))||_{L^2(X)}^2
\\\nonumber&& \leq \int_0^{\tau} \int_X M_{\omega} e^{\omega(\tau - s)}|\chi_B(x) (u_n(x,s) - u(x,s))|^2 dx ds
\\\nonumber&& \leq C(M,\tau) \int_0^{\tau} \int_X |\chi_B(x) (u_n(x,s) - u(x,s))|^2 dx ds
\\\nonumber&& = C(M,\tau) ||(u_n(x,s) - u(x,s))||_{L^2([0,\tau]:L^2(B))}^2 \rightarrow 0.
\end{eqnarray}
This shows that the range of the controllability operator $\mathcal{B}^{\tau}$ is closed. Hence we have exact controllability in
$L^2(\mathcal{R}^{\tau})$.
\end{pf}
The objective of this paper is to provide a solution to the optimal actuator placement problem and hence the optimal location of the set $B$. This motivates us to consider the following definition of controllability gramian parameterized over set $B$.
\begin{defn}
The finite time controllability gramian ${\cal C}_B^{\tau}:L^2(X) \rightarrow L^2(X)$ for the PDE (\ref{transportcontrol}) is given by the following:
\begin{eqnarray}
\label{controllabilitygramian}
{\cal C}_B^{\tau} z = \mathbf{\mathcal{B}}^{\tau} \mathbf{\mathcal{B}}^{\tau *} z
=  \int_0^{\tau} \mathbb{P}_{(\tau - s)} (\chi_B(x) \mathbb{U}_{(\tau - s)} z(x))ds.
\end{eqnarray}
Furthermore, we have the following definition for the induced two norm of the operator ${\cal C}_B^{\tau}$:
\begin{eqnarray}
\nonumber&&||{\cal C}_B^{\tau}||_2^2 = \max_{z\in L^2(X), s.t. \parallel z \parallel_{L^2(X)}=1}\left<{ \cal C}_B^{\tau} z,z\right>_{L^2(X)}.
\end{eqnarray}
\end{defn}
\begin{thm}\label{gramian_theorem}The controllability gramian ${\cal C}^\tau_B : L^2(X)\to L^2(X)$ can be written as a multiplication operator as follows:
\begin{eqnarray}
({\cal C}^\tau_B z)(x)=\left(\int_0^\tau \mathbb{P}_t \chi_B(x)dt\right) z(x)\label{gramian_formula}
\end{eqnarray}
\end{thm}
\begin{pf}
\begin{eqnarray}
\nonumber
&&{\cal C}_B^{\tau} z = \int_0^{\tau} \mathbb{P}_{(\tau - s)} (\chi_B(x) \mathbb{U}_{(\tau - s)} z(x))ds
\\\nonumber&&=  \int_0^{\tau} \mathbb{P}_{s} (\chi_B(x) \mathbb{U}_{s} z(x))ds = \int_0^{\tau} \mathbb{P}_{s} (\chi_B(x) z(\phi_s(x)))ds
\\\nonumber&& =\int_0^{\tau}  \chi_B(\phi_{-s}(x) z(x) \left|\frac{\partial \phi_s(x)}{\partial x}\right|^{-1}ds
= \left[\int_0^{\tau}  (\mathbb{P}_{s}\chi_B(x)) ds\right] z(x).
\end{eqnarray}
\end{pf}
The explicit formula for the controllability gramian from Eq. (\ref{gramian_formula}) in terms of multiplication operator can be used to provide an analytical expression for the minimum energy control input.
\begin{claim} \label{claimcontrollabilitygramian}$\rho_B^\tau(x):=\int_0^\tau \mathbb{P}_t \chi_B(x)dt$ is strictly positive on ${\cal R}^\tau =\cup_{t=0}^{\tau} \phi_t(B)$ and hence ${\cal C}_B^\tau$ is invertible on $\mathcal{R}^\tau$ with the inverse given by
\begin{equation} \label{inverseofcontrollabilitygramian}({\cal C}^{\tau}_B)^{-1} z=\frac{z}{\rho_B^\tau(x)},\;\;\;\;\forall z\in L^2({\cal R}^\tau).\end{equation}
\end{claim}
\begin{pf}
Since $m(B) > 0$, and $B$ evolves into $\phi_{\tau}(B)$ in time $\tau$,
for every $x \in \mathcal{R}^{\tau}$, there exist times $0 \leq t_1(x) <  t_2(x) \leq \tau$ such that
$x \in \phi_t(B) ~\forall t \in [t_1(x),t_2(x)]$. Hence, by the positivity of $\mathbb{P}_t$ we have that
$\mathbb{P}_{t}(\chi_B(x)) > 0 ~\forall t \in [t_1(x),t_2(x)] \subseteq [0,\tau]$. Hence we have the following:
\begin{eqnarray}
\nonumber \rho_B^\tau(x)=\int_0^\tau \mathbb{P}_t \chi_B(x)dt \geq \int_{t_1(x)}^{t_2(x)} \mathbb{P}_t \chi_B(x)dt > 0 ~\forall x \in \mathcal{R}^{\tau}.
\end{eqnarray}
This proves the claim.
\end{pf}
\begin{thm}\label{theorem_minimum_energy}
Let $\rho_{\tau}(x)$ and $\rho_0(x)$ be the elements of $L^2({\cal R}^\tau)$, then the minimum energy control input that is required to steer the system from initial state $\rho_0(x)$ to final state $\rho_{\tau}(x)$ is given by following formula
\begin{eqnarray}
\nonumber&& u_{\mbox{{\it opt}}}(x,s)= \mathcal{B}^{\tau *}(\mathcal{C}_{B}^{\tau})^{-1}(\rho_{\tau}(x)-\mathbb{P}_{\tau} \rho_0(x))
\\\label{minimumenergycontrolinput}&& =  \chi_B(x)\mathbb{U}_{\tau-s} \left(\frac{\rho_{\tau}(x)-\mathbb{P}_{\tau} \rho_0(x)}{\rho_B^{\tau}(x)}\right).
\end{eqnarray}
The minimum energy required is given by
\begin{eqnarray}\label{minimumenergycontrolinputnorm}&&||u_{\mbox{{\it opt}}}||^2
\\\nonumber&&= \left<(\rho_{\tau}(x)-\mathbb{P}_{\tau} \rho_0(x)),(\mathcal{C}_B^{\tau})^{-1} (\rho_{\tau}(x)-\mathbb{P}_{\tau} \rho_0(x))\right>_{L^2(\mathcal{R}^{\tau})}
\\\nonumber&&=\left|\left|\frac{(\rho_{\tau}(x)-\mathbb{P}_{\tau} \rho_0(x))}{\rho_B^{\tau}(x)}\right|\right|^2_{L^2(\mathcal{R}^{\tau})}.
\end{eqnarray}
\end{thm}
\begin{pf}
First, we note that controlling the initial state $\rho_0(x)$ to $\rho_{\tau}(x)$ is equivalent to reaching the final state $(\rho_{\tau}(x)-\mathbb{P}_{\tau}\rho_0(x))$
from the zero initial state i.e. $\rho_0(x) \equiv 0$. Hence, equivalently, we prove that $\hat{u}_{\mbox{{\it opt}}}(x,s)= \mathcal{B}^{\tau *}(\mathcal{C}_{B}^{\tau})^{-1}(\rho_{\tau}(x))$
is the control input with minimum norm that reaches $\rho_{\tau}(x)$ in time $\tau$. This, along with an explicit calculation of $\mathcal{B}^{\tau *}(\mathcal{C}_{B}^{\tau})^{-1}(\rho_{\tau}(x))$
will prove the Theorem.
Next, we consider the following set of admissible control inputs:
\[\mathcal{U} = \{u(x,t) \in L^2([0,\tau]:L^2(B)): \mathcal{B}^{\tau} u = \rho_{\tau}\}.\]
We have the following:
\begin{eqnarray}
\nonumber \mathcal{B}^{\tau} \hat{u}_{\mbox{{\it opt}}} = \mathcal{B}^{\tau} \mathcal{B}^{\tau *}(\mathcal{C}_{B}^{\tau})^{-1}\rho_{\tau} =
\mathcal{B}^{\tau} \mathcal{B}^{\tau *}(\mathcal{B}^{\tau} \mathcal{B}^{\tau *})^{-1} \rho_{\tau} = \rho_{\tau}.
\end{eqnarray}
Hence, we have that
$\hat{u}_{\mbox{{\it opt}}}(x,s)= \mathcal{B}^{\tau *}(\mathcal{C}_{B}^{\tau})^{-1}(\rho_{\tau}(x)) \in \mathcal{U}$.
Next, we define the following operator on $L^2([0,\tau]:L^2(B))$ $P^{\tau} = \mathcal{B}^{\tau *}(\mathcal{C}_{B}^{\tau})^{-1}\mathcal{B}^{\tau}$.
We observe the following:
\begin{eqnarray}
(P^{\tau})^2 &=& \mathcal{B}^{\tau *}(\mathcal{C}_{B}^{\tau})^{-1}\mathcal{B}^{\tau}\mathcal{B}^{\tau *}(\mathcal{C}_{B}^{\tau})^{-1}\mathcal{B}^{\tau}=\mathcal{B}^{\tau *}(\mathcal{C}_{B}^{\tau})^{-1}\mathcal{B}^{\tau}\nonumber
\\&=& P^{\tau},(P^{\tau})^* = (\mathcal{B}^{\tau *}(\mathcal{C}_{B}^{\tau})^{-1}\mathcal{B}^{\tau})^* = P^{\tau}.
\end{eqnarray}
Hence, the operator $P^{\tau}$ is a projection operator on the space $L^2([0,\tau]:L^2(B))$. Then, we have the following from Bessel's inequality:
\[||u||^2 = ||(P^{\tau})u||^2 + ||(I - P^{\tau})u||^2 \geq ||(P^{\tau})u||^2,\]
where the norm is on the space $L^2([0,\tau]:L^2(B))$. Now, let $u \in \mathcal{U}$ be arbitrary. This means
$\mathcal{B}^{\tau} u = \rho_{\tau}$. Applying $\mathcal{B}^{\tau *} (\mathcal{C}_B^{\tau})^{-1}$ on both sides, we get the following:
\[P^{\tau} u = \mathcal{B}^{\tau *} (\mathcal{C}_B^{\tau})^{-1} \mathcal{B}^{\tau} u = \mathcal{B}^{\tau *} (\mathcal{C}_B^{\tau})^{-1} \rho_{\tau} = \hat{u}_{\mbox{{\it opt}}}.\]
Hence, Bessel's inequality above gives $||u||^2 \geq ||\hat{u}_{\mbox{{\it opt}}}||^2$. Next, (\ref{minimumenergycontrolinput}) and
(\ref{minimumenergycontrolinputnorm}) can be easily shown by an explicit calculation using (\ref{adjointofcontrollabilityoperator})
and (\ref{inverseofcontrollabilitygramian}).
\end{pf}
Based on the formula for the controllability gramian, we propose the following criteria for the selection of optimal actuator location and hence the set $B^*$. \\
{\it Actuator placement criteria}
\begin{enumerate}
\item Maximizing the support of the controllability gramian operator i.e.,
\begin{eqnarray}
B^{*}&=&\operatorname*{arg\max}_{B\subset X}\;\;supp \left(\int_{0}^{\tau}\mathbb{P}_t \chi_B(x)dt \right)
%\nonumber\\
%&=&\operatorname*{arg\max}_{B\subset X}\;\;supp \left(\int_{0}^{\tau} \chi_{\phi_t(B)}(x) \left|\frac{\partial \phi_t(x)}{\partial x}\right|^{-1}dt \right)\nonumber
\end{eqnarray}
\item If the support of controllability gramian is maximized or if more than one choice of set $A$ leads to the same support then the decision can be made based on maximizing the $2$-norm of the support i.e.,

\[B^{*}=\operatorname*{arg\,max}_{B\subset X} \parallel \int_0^\tau \mathbb{P}_t \chi_B(x) dt \parallel_{L^2(X)}.\]

\end{enumerate}
%We now consider two different type of vector fields $f(x)$, for which we will compute the infinite time controllability gramian.
Using the result of Theorem \ref{pde_controllability}, it follows that criterion 1 maximizes the controllability in the space $X$, so that the control action in a small region $B\subset X$ will have an impact over larger portion of the state space. Furthermore, it follows from the explicit formula for the minimum energy control (\ref{minimumenergycontrolinput}) from Theorem \ref{theorem_minimum_energy} that if the the actuator selection is made based on criteria 2 then the amount of control effort is minimized.

\subsection{Observability gramian}

For the construction of observability gramian, we consider the advection partial differential equation with a single output measurement as follows:
\begin{eqnarray}
\frac{\partial \rho}{\partial t}&=&\nabla \cdot (f \rho), \;\;\;\rho|_{\Gamma_i}=0,\;\;\;\rho(x,0)=\rho_0(x)\nonumber\\
y(x,t)&=&\chi_A(x)\rho(x,t)\label{system}
\end{eqnarray}
The observability operator ${\cal A}^{\tau}: L^2(X)\to L^2([0,\tau]:L^2(A))$ for (\ref{system}) is defined as follows:
\[({\cal A}^\tau z)(x,s)=\chi_A(x) (\mathbb{P}_s z)(x). \]
The adjoint to the observability operator ${\cal A^{\tau *}} : L^2([0,\tau]:L^2(A)) \to L^2(X)$ can be written as follows:
\[({\cal A}^{\tau *} w)(x)=\int_{0}^{\tau} (\mathbb{U}_{s} \chi_A(x) w(x,s))ds.\]
\begin{defn}[Observability gramian] The finite time observability gramian ${\cal O}_A^{\tau}: L^2(X)\to L^2(X)$ for the PDE (\ref{system}) is given by the following formula
\begin{eqnarray}
({\cal O}_A^{\tau} z)(x)=({\cal A}^{\tau *} {\cal A}^\tau  z)(x)=\int_0^{\tau}(\mathbb{U}_s \chi_A(x) \mathbb{P}_s z(x))ds. \label{obsver_gramian_eqn}
\end{eqnarray}
\end{defn}
The  counterpart of Theorems  (\ref{pde_controllability}) and (\ref{theorem_minimum_energy}) can be proved for the observability of system (\ref{system}) using a duality argument. The theorem on observability gramian similar to Theorem (\ref{gramian_theorem}) can be stated as follows:
\begin{thm}\label{observabilitygramiantheorem} The observability gramian for (\ref{system}) can be written as a multiplication operator as follows:
\begin{eqnarray}
({\cal O}_A^{\tau} z)(x)=\left[\int_0^\tau (\mathbb{U}_s \chi_A(x)) ds\right] z(x).
\end{eqnarray}
\end{thm}
\begin{pf}
The proof follows along the lines of proof of Theorem (\ref{gramian_theorem}).
\end{pf}
Following criteria can be used for the optimal location of sensor:

\noindent{\it Sensor placement criteria}\\
The finite time observability gramian can be used to decide the criteria for the optimal location of the sensor.
\begin{enumerate}
\item Maximizing the support of observability gramian operator
\[A^{*}=\operatorname*{arg\,max}_{A\subset X} supp \left(  \int_0^\tau \mathbb{U}_t \chi_A(x) dt\right ).\]

\item If the support of observability gramian is maximized or if more than one choice of set $B$ leads to the same support then the decision can be made based on maximizing the $2$-norm of the support i.e.,

\[A^{*}=\operatorname*{arg\,max}_{A\subset X} \parallel \int_0^\tau \mathbb{U}_t \chi_A(x) dt \parallel_{L^2(X)}\]

\end{enumerate}

\section{Advective vector field and gramian}
\label{Infinitetimegramians}
In this section, we provide an interpretation for the optimal actuator and sensor location problem in terms of the flow of the advection vector field. In particular, we show that the (almost everywhere uniform) stability property of the vector field plays an important role in  deciding the location of actuators and sensors.

\subsection{Infinite time Controllability gramian}
We show that the infinite time controllability gramian  can be computed for vector
fields that are stable in the almost everywhere uniform sense.
We now define the notion of almost everywhere uniform stability for a nonlinear system.
\begin{defn}[Almost everywhere uniform stable]\label{assume_stability} Let $x_0=0$ be the equilibrium point of $\dot x=f(x)$ and $B_{\delta}$ be a $\delta$ neighborhood of $x_0=0$.
The equilibrium point $x_0=0$ is said to be almost everywhere uniform stable if  for every given $\epsilon>0$ there exists a $T(\epsilon)$ such that
\[\int_T^{\infty}m(A_t)dt<\epsilon,;\;\;\;\;A_t=\{x\in X: \phi_t(x)\in A\},\]
for all measurable sets $A\subset X\setminus B_{\delta}$ and where $m$ is the Lebesgue measure.
\end{defn}
The notion of almost everywhere stability is extensively studied in \cite{Rantzer01} \cite{Vaidya_TAC}. Furthermore, a PDE based approach is also provided for the verification of almost everywhere stability in \cite{stabilitypde}.  We have the following theorem regarding the infinite time controllability gramian for vector fields that are almost everywhere uniformly stable:
\begin{thm}\label{infinitetimecontrollabilitygramian} For vector fields that are stable in the almost everywhere uniform sense, we have
\begin{eqnarray}
({\cal C}^\infty_Bz)(x)=\int_{0}^\infty \mathbb{P}_t\chi_B(x)dt z(x)=\rho_B(x)z(x),\label{infinite_controllability}
\end{eqnarray}
where $\rho_B(x)$ is the positive solution  of the following PDE
\begin{equation} \label{pde} \nabla \cdot (f (x)\rho_B(x))=\chi_B(x);\rho|_{\Gamma_i} = 0.\end{equation}
\end{thm}
\begin{pf}
In \cite{stabilitypde}, it was shown that $\int_{0}^\infty \mathbb{P}_t\chi_B(x)dt$ solves (\ref{pde}) if $x_0=0$ is stable in the
almost everywhere uniform sense. This proves the Theorem.
\end{pf}
The integral $\int_X {\cal C}^\infty_B z(x) dx$ for the special case where $z(x)=\chi_A(x)$, the indicator function for the set $A$, has the interesting interpretation of residence time, which is defined as follows:
\begin{defn} For an almost everywhere uniform stable vector field, consider any two measurable subsets $A$ and $B$ of $X\setminus B_{\delta}$, then the residence time of set B in set $A$ is defined as the amount of time system trajectories starting from set $B$ will spend in set $A$ before entering the $\delta$ neighborhood of the equilibrium point $x=0$. We denote this time by $T_B^A$.
\end{defn}
In \cite{Umeshresidencetimes}, the following was shown for a discrete time system:
\begin{equation} \label{residencetimeequation} T_B^A=\int_A\int_0^{\infty} \mathbb{P}_t \chi_B(x)dt dx = \int_A \rho_B(x) dx. \end{equation}
The proof of (\ref{residencetimeequation}) for a continuous-time case will follow along the lines of proof in \cite{Umeshresidencetimes}.
\begin{thm}\label{residence_time}
The residence time $T_B^A$ for an almost everywhere uniformly stable vector field $f(x)$ is given by following formula
\[T_B^A=\int_X {\cal C}^\infty_B \chi_A(x) dx.\]
\end{thm}
\begin{pf}
We have the following calculation using the formula from Theorem \ref{infinitetimecontrollabilitygramian} and
Eq. (\ref{residencetimeequation}):
\begin{eqnarray}
\nonumber&& ({\cal C}^\infty_B \chi_A(x))=\int_{0}^\infty \mathbb{P}_t\chi_B(x)dt \chi_A(x) = \rho_B(x) \chi_A(x)
\\\nonumber&& \Rightarrow \int_X {\cal C}_B^{\infty}\chi_A(x)  dx = \int_X \rho_B(x) \chi_A(x) dx
= \int_A \rho_B(x) dx = T_B^A.
\end{eqnarray}
\end{pf}
\subsection{Infinite time observability gramian}
The infinite time observability gramian is defined under the assumption that the vector field $f(x)$ is globally asymptotically stable.
First, we have the following Theorem that characterizes global asymptotic stability:
\begin{thm}
\label{steadystatetransporttheorem}
Let $B_{\delta}$ be a $\delta$ neigborhood of $x=0$.
Let $v(x) \in C^1(X \setminus \bar{B}_{\delta})$ denote the solution of the following steady state transport equation:
\begin{equation}
\label{steadystatetransportequation}
\mathcal{A}_Kv = f \cdot \bigtriangledown v = - v_0(x); v|_{\partial \bar{B}_{\delta}} = 0,
\end{equation}
where $v_0(x)$ satisfies
\begin{equation} \label{v0condition} v_0(x)=0~\forall x \in \bar{B}_{\delta}. \end{equation}
Then $x=0$ is globally asymptotically stable for (\ref{ode}) if and only if there exists a
positive solution $v(x) \in C^1(X/ \bar{B}_{\delta})$ for (\ref{steadystatetransportequation})
for all $v_0(x)>0 \in C^1(X/ \bar{B}_{\delta})$ satisfying (\ref{v0condition}).
\end{thm}
\begin{pf}
We prove necessity first. Let us assume that $x=0$ is globally asymptotically stable. We construct a positive solution for
(\ref{steadystatetransportequation}) as follows:
\begin{equation}
\label{steadystatetransportsolution}
v(x) = \int_0^{\infty} v_0(\phi_t(x)) dt.
\end{equation}
For a given arbitrary $x \in X/\bar{B}_{\delta}$, there exists a $t^+(x) \in [0,\infty)$ such that $\phi_{t^+(x)}(x) \in \partial B_{\delta}$.
Hence, we have that $v_0(\phi_t(x))=0 ~\forall t \geq t^+(x)$. In particular, this means that $\int_0^{\infty} v_0(\phi_t(x)) dt < \infty
~\forall x \in X/\bar{B}_{\delta}$. We also have that $0<v(x) \in C^1(X/\bar{B}_{\delta})$ by virtue of the regularity of
$v_0(x)$. We show that (\ref{steadystatetransportsolution}) solves (\ref{steadystatetransportequation}). Let
$v_N(x) = \int_0^N v_0(\phi_t(x)) dt$. Then, we have the following:
\begin{eqnarray}
\nonumber&& \mathcal{A}_K v_N(x) = \int_0^N \mathcal{A}_K v_0(\phi_t(x)) dt =
\\\nonumber&&\int_0^N \frac{d}{dt} \mathbb{U}_t v_0(x) dt = \mathbb{U}_N v_0(x) - v_0(x).
\end{eqnarray}
Global stability of $x=0$ implies that $\displaystyle \lim_{N \rightarrow \infty}\mathbb{U}_N v_0(x) = \lim_{N \rightarrow \infty} v_0(\phi_N(x)) = 0$ and hence
$\displaystyle \lim_{t \rightarrow \infty} \mathcal{A}_K v_N(x)$ exists. Also, by the Hille-Yosida
semigroup generation theorem, we have that the generator $\mathcal{A}_K$ is a closed operator. Hence,
we have the following:
\begin{eqnarray}
\nonumber f \cdot \bigtriangledown  v =
\mathcal{A}_K v(x) = \int_0^{\infty} \mathcal{A}_K \mathbb{U}_t v_0(x) = \int_0^{\infty} \frac{d}{dt} \mathbb{U}_t v_0(x) = -v_0(x).
\end{eqnarray}
The boundary condition $v|_{\partial \bar{B}_{\delta}} = 0$ is satisfied by (\ref{steadystatetransportsolution}) automatically.
To prove sufficiency let us assume that there exists a solution $0<v(x) \in C^1(X/\bar{B}_{\delta})$ that solves (\ref{steadystatetransportequation}).
Then, we have the following equation along the characteristic curves which are solutions of (\ref{ode}):
\begin{eqnarray}
\label{characteristicequation}&&
\frac{d}{d\tau} v(\phi_{\tau}(x)) = -v_0(\phi_{\tau}(x)) \Rightarrow v(\phi_t(x)) - v(x)
\\\nonumber&&= -\int_0^t v_0(\phi_{\tau}(x)) d\tau.
\end{eqnarray}
Rewriting (\ref{characteristicequation}), we have $v(\phi_t(x)) + \int_0^t v_0(\phi_{\tau}(x)) d\tau = v(x)$
\vspace{-0.05in}
\begin{eqnarray}%
\nonumber&& \Rightarrow \int_0^t v_0(\phi_{\tau}(x)) d \tau \leq v(x) ~\forall x \in X/\bar{B}_{\delta}, t>0
\\\label{eqn1}&& \Rightarrow ||\int_0^{\infty} v_0(\phi_{\tau}(x)) d \tau ||_{L^{\infty}(X/\bar{B}_{\delta})}
\leq ||v(x)||_{L^{\infty}(X/\bar{B}_{\delta})} < \infty.
\end{eqnarray}
To the contrary, let us assume that $x=0$ is not globally asymptotically stable. Then,
by virtue of the attractor property of $x=0$, there exists a point $x_0 \in X/ \bar{B}_{\delta}$ such that $\omega(x_0) \neq \{0\}$. This means
that $\phi_t(x_0) \in X / \bar{B}_{\delta} ~\forall t>0$, for some $\delta > 0$. Then, the set
$D=\overline{\cup_{t=0}^{\infty} \phi_{t}(x_0)}$ is a compact subset of $X/\bar{B}_{\delta}$. Since $v_0(x) > 0 ~\forall x \in X/ \bar{B}_{\delta}$, we have that
$v_0(x) > \epsilon > 0 ~\forall x \in D$ for some positive $\epsilon$ by continuity of $v_0(x)$. Hence we have the following:
\begin{equation}
\label{contradiction2} \int_0^{\infty} v_0(\phi_{\tau}(x_0)) d \tau > \int_0^{\infty} \epsilon d \tau = \infty,
\end{equation}
contradicting (\ref{eqn1}). This proves the Theorem.
\end{pf}
If $x=0$ is globally asymptotically stable, then $\Gamma_o \supseteq \partial \bar{B}_{\delta}$.
Hence, by using a standard density argument of $C^1(X \setminus \bar{B}_{\delta})$ in $L^2(X \setminus B_{\delta})$, and
using trace operator theory \cite{Evans} for point values of $H^1$ functions, we can show the following Theorem:
\begin{thm}
\label{L2theorem}
Let $v(x) \in \mathcal{D}(\mathcal{A}_K) \cap L^2(X \setminus \bar{B}_{\delta})$ denote the solution of the following steady state transport equation:
\begin{equation}
\label{steadystatetransportequationL2}
\mathcal{A}_K v = f \cdot \bigtriangledown v = - v_0(x); v|_{\Gamma_o} = 0,
\end{equation}
Then $x=0$ is globally asymptotically stable for (\ref{ode}) if and only if there exists a
positive solution $v(x) \in \mathcal{D}(\mathcal{A}_K) \cap L^2(X/ \bar{B}_{\delta})$ for (\ref{steadystatetransportequation})
for all $0 < v_0(x) \in \mathcal{D}(\mathcal{A}_K) \cap L^2(X/ \bar{B}_{\delta})$.
\end{thm}
\begin{thm}Let $x=0$ be a globally stable equilibrium point for $\dot x=f(x)$, then the infinite time observability gramian is well defined and we have
\begin{eqnarray}
\label{infinitetimeobservabilitygramianequation}
({\cal O}^\infty_A z)(x) = \left[\int_0^\infty (\mathbb{U}_t \chi_A(x))dt \right] z(x)=V(x)z(x),\label{infinite_observability}
\end{eqnarray}
where $V(x)$ is the positive solution of following steady state partial differential equation:
\[\mathcal{A}_K v = f \cdot \bigtriangledown v = - \chi_A(x); v|_{\Gamma_o} = 0.\]
\end{thm}
\begin{pf}
For a given $A$, if we choose $\delta > 0$ such that $A \subset X \setminus B_{\delta}$, then we automatically have that
$\chi_A(x) = 0 ~\forall x \in \bar{B}_{\delta}$. Hence, global stability implies the existence of a positive solution
$V(x) = \int_0^\infty (\mathbb{U}_t \chi_A(x))dt \in \mathcal{D}(\mathcal{A}_K) \cap L^2(X \setminus B_{\delta})$ from Theorems \ref{steadystatetransporttheorem} and(\ref{L2theorem}).
This shows that $V(x)$ is well defined. Finally, the formula for the infinite time observability gramian (\ref{infinitetimeobservabilitygramianequation}) is obtained by letting
$\tau \rightarrow \infty$ in Theorem \ref{observabilitygramiantheorem}. This proves the Theorem.
\end{pf}
%\subsection{Fluid flow vector field}

%\input{complex_field_umesh}
\section{Simulation}\label{simulation}
\begin{figure}[h*]
\begin{center}
\mbox{
\hspace{-0.2in}
\subfigure[]{\scalebox{1}{\includegraphics[width=2.3 in, height=2.0 in]{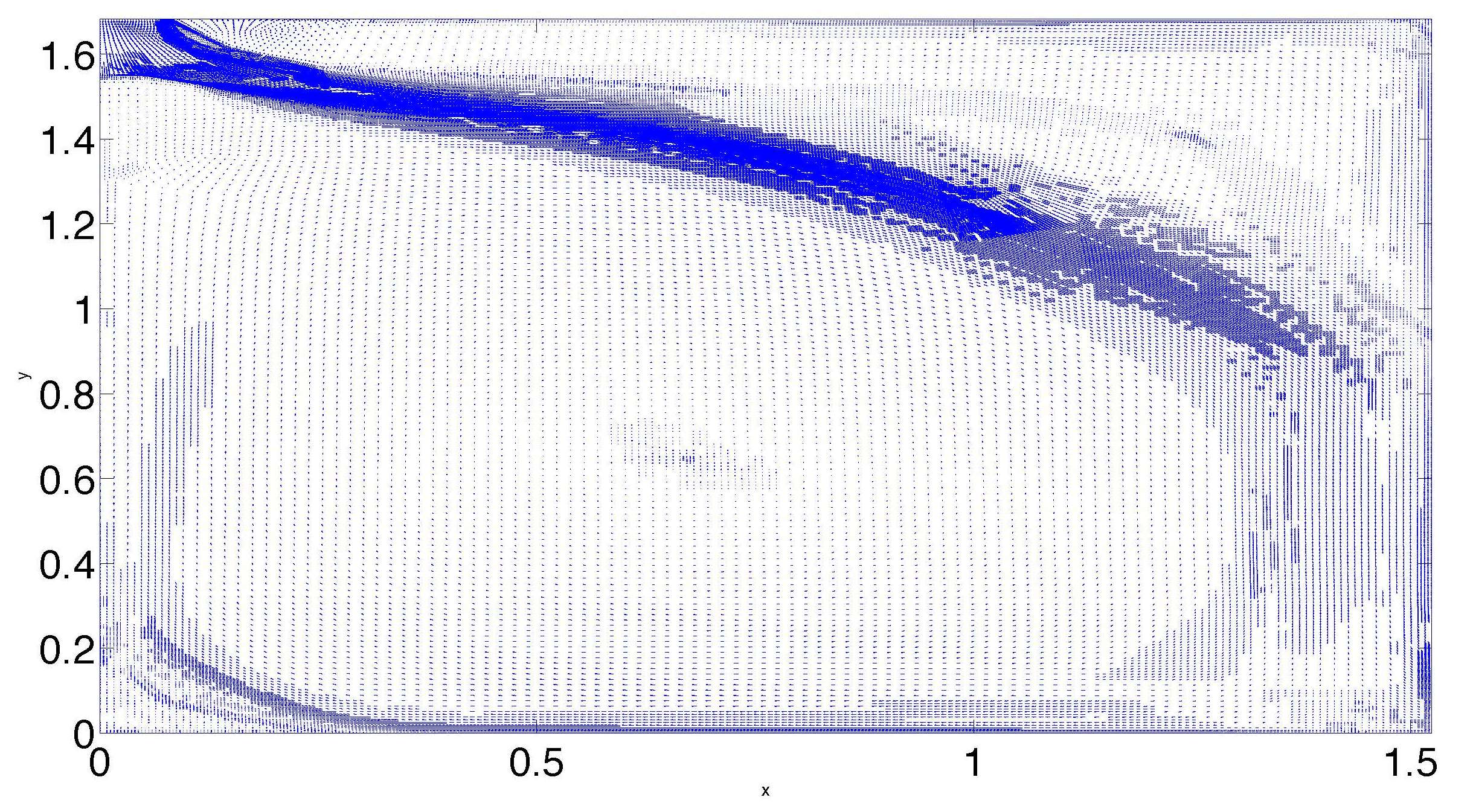}}}
%\subfigure[]{\scalebox{.8}{\includegraphics[width=3.3 in,height=3.0 in]{vf.jpg}}}
\hspace{-0.0in}
\subfigure[]{\scalebox{.35}{\includegraphics[width=6.3 in, height=6.0 in]{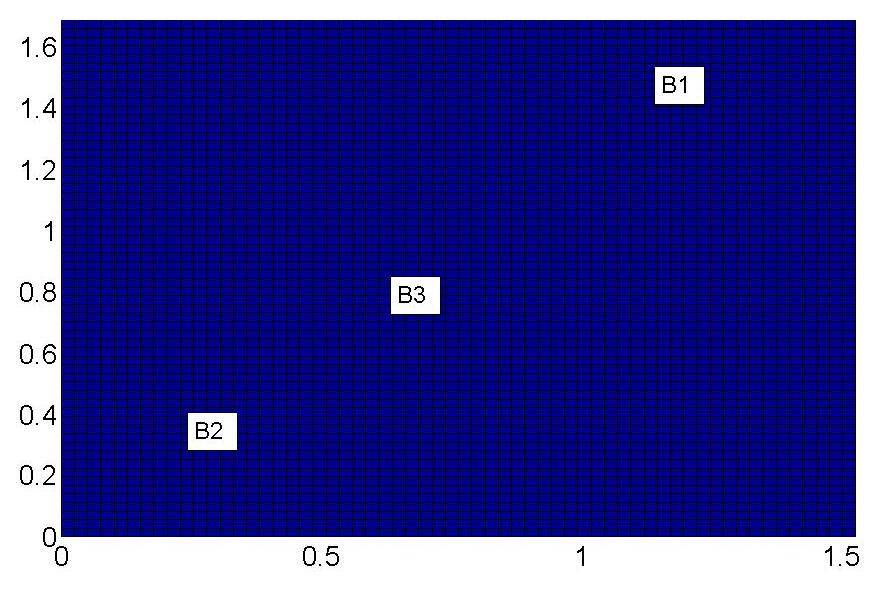}}} }
%\subfigure[]{\scalebox{.26}{\includegraphics[width=7 in, height=6.7 in]{sets.jpg}}} }
\vspace{-0.1in}
\caption{(a)  Velocity field; b) Actuator locations on sets $B_1, B_2$ and $B_3$.}
\label{vec}
\end{center}
\end{figure}

In this section, we present simulation results on the computation of finite time gramians.  The purpose of the simulation section is to demonstrate the applicability of the developed theoretical results in this paper. Detailed simulation results based on the developed theoretical results will be the topic of our future publication. The vector field that we use for the purpose of simulation is the average velocity field obtained from a detailed finite element-based simulation of Navier Stokes equation. For the purpose of simulation, we only employ a two dimensional slice of the three dimensional velocity field as shown in Fig. \ref{vec}a. The dimensions of the room are as follows: $0\leq x\leq 1.52 m$ and $0\leq y\leq 1.68 m$. The order of magnitude for the velocity field is $O(1)$. The Reynolds number of the flow is $Re=76725$ and the Prandtl number $Pr=0.729$. This makes $\frac{1}{Pr Re}\approx O(10^{-5})$, and hence the zero diffusion constant assumption (Assumption \ref{assume_zero_diffusion}) made in this paper is justified. The Reynolds number for the flow rate is in turbulent range. The $k-\epsilon$ model, which is Reynolds Average Navier-Stokes (RANS) model \cite{k-eplsion_model} is used to obtain the velocity field as shown in Fig. \ref{vec}. A commercial CFD software Fluent was used to solve  the coupled set of governing equations for pressure, temperature, turbulent kinetic energy, turbulent dissipation and velocity. No slip boundary condition was applied at all the walls.

For the purposes of computation, we employ set oriented numerical methods for the approximation of P-F semigroup $\mathbb{P}_t$ \cite{Dellnitz00}. We
divide the state space into finitely many square partitions denoted by
$\{D_i\}_{i=1}^N$. The set $D_i$'s are chosen such that $D_i\cap D_j=\emptyset$ for $i\neq j$ and $X=\cup_{i=1}^N D_i$. The finite dimensional matrix approximation of the P-F operator is obtained using the following formula \cite{Dellnitz00}:
\[[P]_{ij}=\frac{m(\phi_{\delta t}(D_i)\cap D_j)}{m(D_i)}\]
where $m$ is the Lebesgue measure and $\delta t$ is the discretization time step and $\phi_t$ is the solution of vector field shown in Fig. \ref{vec}a. Using the adjoint property between the Koopman and P-F semigroup, the finite dimensional approximation of the Koopman semigroup $U$ can be obtained as a transpose of $P$, namely $U=P'$

The computation results for this section are obtained with actuators and sensors located at three different sets $B_1,B_2$, and $B_3$. The locations of these three sets are shown in Fig. \ref{vec}b.
%In this section, we verify the theoretical findings with numerical simulations. The quantity, which is to be computed and which determines controllability for time $\tau$ is $\rho_B^\tau(x):=\int_0^\tau \mathbb {P}_t \chi_B (x)dt$. For simulation purposes, first we approximate the continuous time Perron Frobenius operator by the discrete time one as outlined in [ ] and [ ] and then we use Ulam's approximation method and obtain a finite dimensional approximation for the operator $\mathbb{P}_t$ and obtain a finite dimensional approximation. We equally partition the space $X$  as  $X = {D_1, D_2,......D_N}$. Then take $ \chi_{D_i}(x)$ as the basis function for $\rho(x,t)$ and obtain a finite dimensional approximation. The infinite dimensional PF operator becomes a Markov matrix $P$ such that,

%\[P_{ij} = \frac{m(T^{-1}(D_i \cap D_j) } {m(D_i)} \]

%where $x(k+1) = T (x(k)) $ is discrete time approximation of $\dot{x}(t) = f(x(t))$. Next we present simulation results for two different vector fields,

%\subsection{ Room Temperature Advection}

%
%\begin{figure}[H]    \centering
%{\includegraphics[width=2.50 in,height=2.2 in]{vf.jpg}}
% \caption{ Vector field profile} \label {vec}
%\end{figure}

%\begin{figure}[h*]
%\begin{center}
%\mbox{
%\hspace{-0.05in}
%\subfigure[]{\scalebox{.3}{\includegraphics{1a}}}
%\hspace{-0.05in}
%\subfigure[]{\scalebox{.3}{\includegraphics{1g}}} }
%\caption{(a)  A particular Sensor Position ($B_1$) (b) Gramian after time 100}
%\label{room1}
%\end{center}
%\end{figure}

%\begin{figure}[H]  \centering  \subfigure[]
%{\includegraphics[width=2.50 in]{1a}} \subfigure[]
%{\includegraphics[width=2.50 in]{1g}}
% \caption{ \label {room1}  (a)  A particular Sensor Position ($B_1$) (b) Gramian after time 100}
%\end{figure}

\begin{figure}[h*]
\begin{center}
\mbox{
\hspace{-0.2in}
\subfigure[]{\scalebox{.35}{\includegraphics{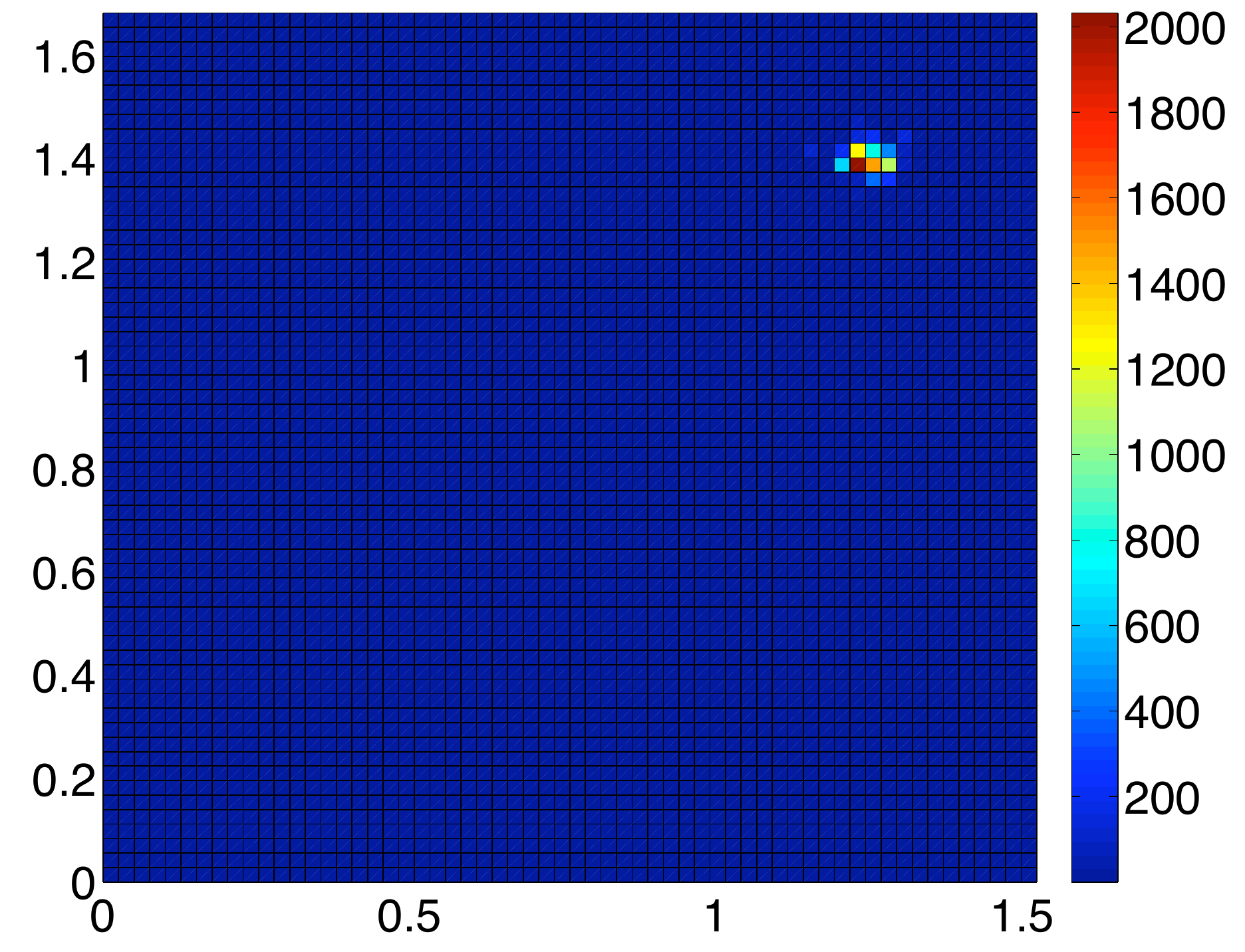}}}
\hspace{-0.2in}
\subfigure[]{\scalebox{.35}{\includegraphics{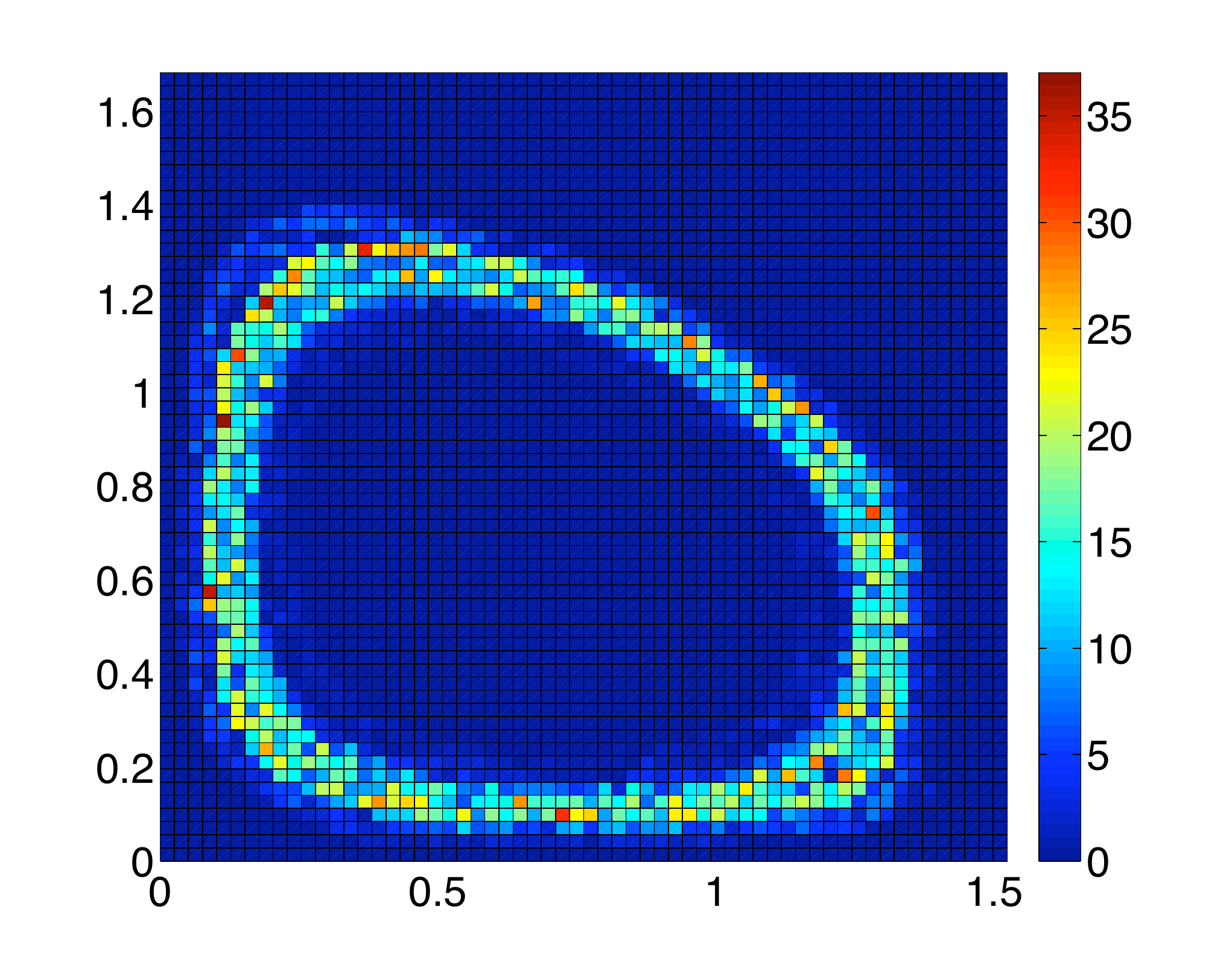}}} }
\vspace{-0.1in}
\caption{ Controllability gramian after 10000 time iterations  for actuator located at set  a)  $B_1$;  b) $B_2$.}
\label{control12}
\end{center}
\end{figure}

%\begin{figure}[H] \centering \subfigure[]
%{\includegraphics[width=2.50 in]{2a}} \subfigure[]
%{\includegraphics[width=2.50 in]{2g}}
% \caption{  \label{room2} (a)  A particular Sensor Position ($B_2$) (b) Gramian after time 100}
%\end{figure}

\begin{figure}[h*]
\begin{center}
\mbox{
\hspace{-0.2in}
\subfigure[]{\scalebox{.35}{\includegraphics{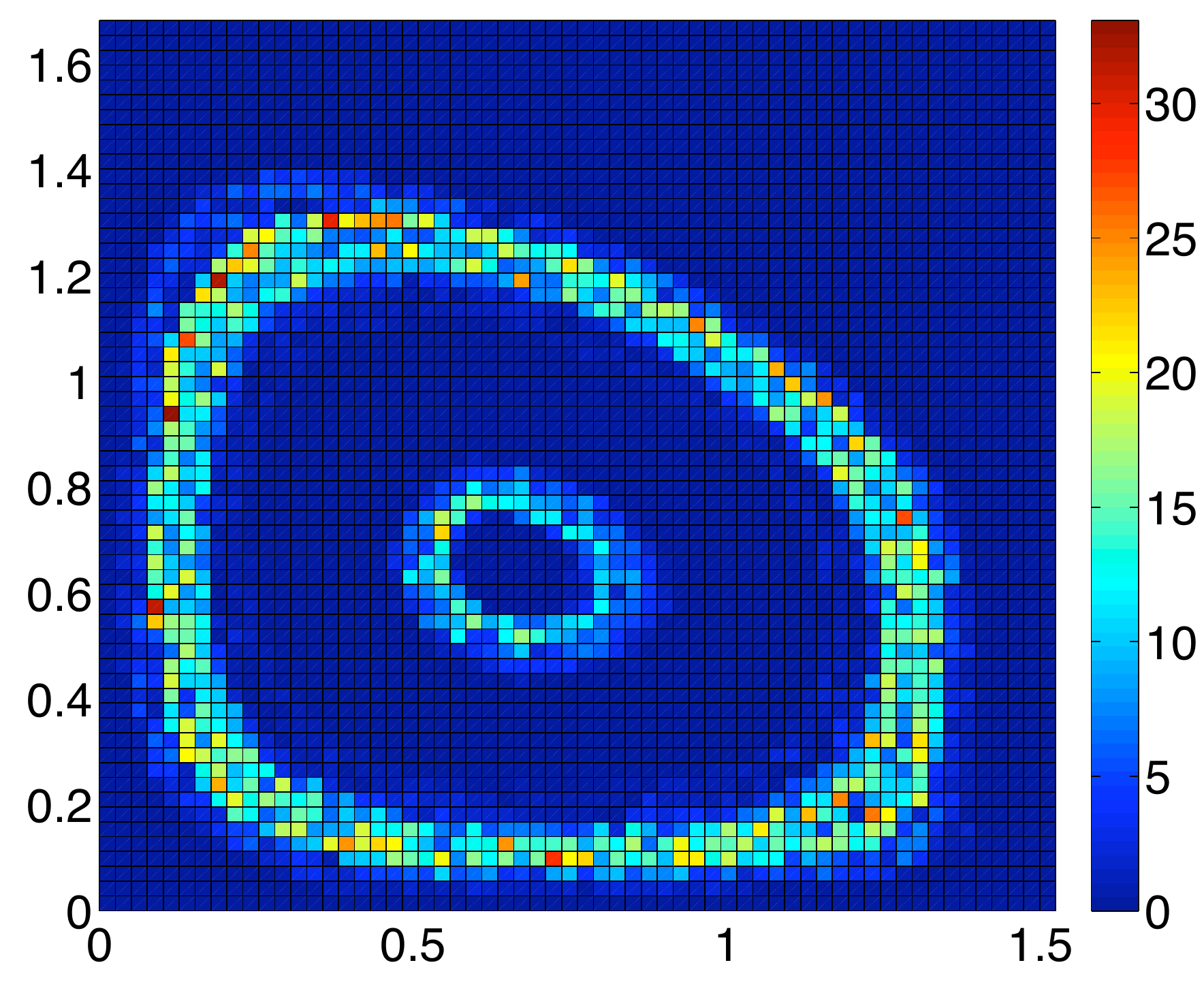}}}
\hspace{-0.2in}
\subfigure[]{\scalebox{.35}{\includegraphics{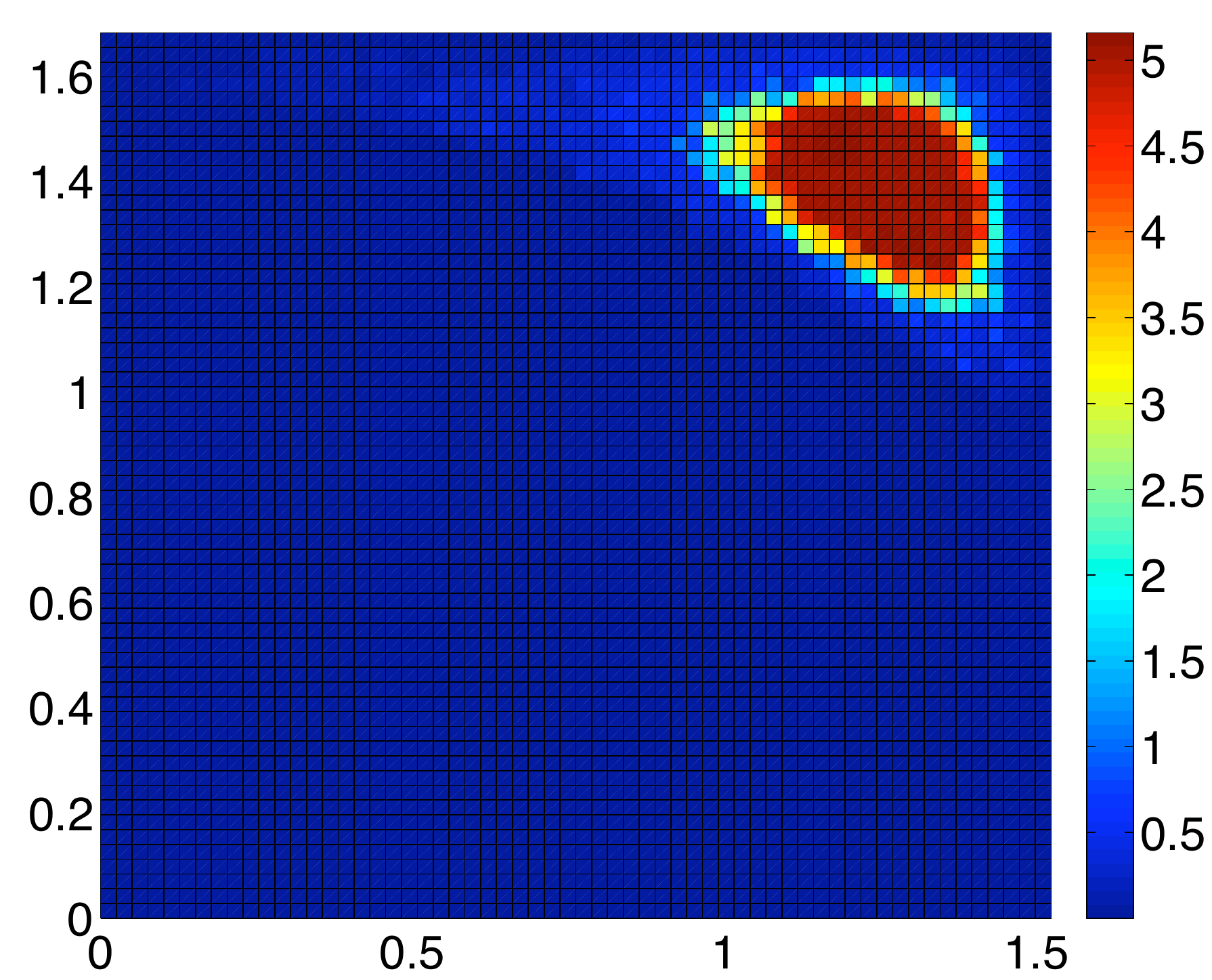}}} }
\vspace{-0.1in}
\caption{a) Controllability gramian after 10000 time iterations  for actuator located at set  $B_3$; b)  Observability gramian after 1000 iterations with sensor location at $B_1$}
\label{control3}
\end{center}
\end{figure}

\begin{figure}[h*]
\begin{center}
\mbox{
\hspace{-0.2in}
\subfigure[]{\scalebox{.35}{\includegraphics{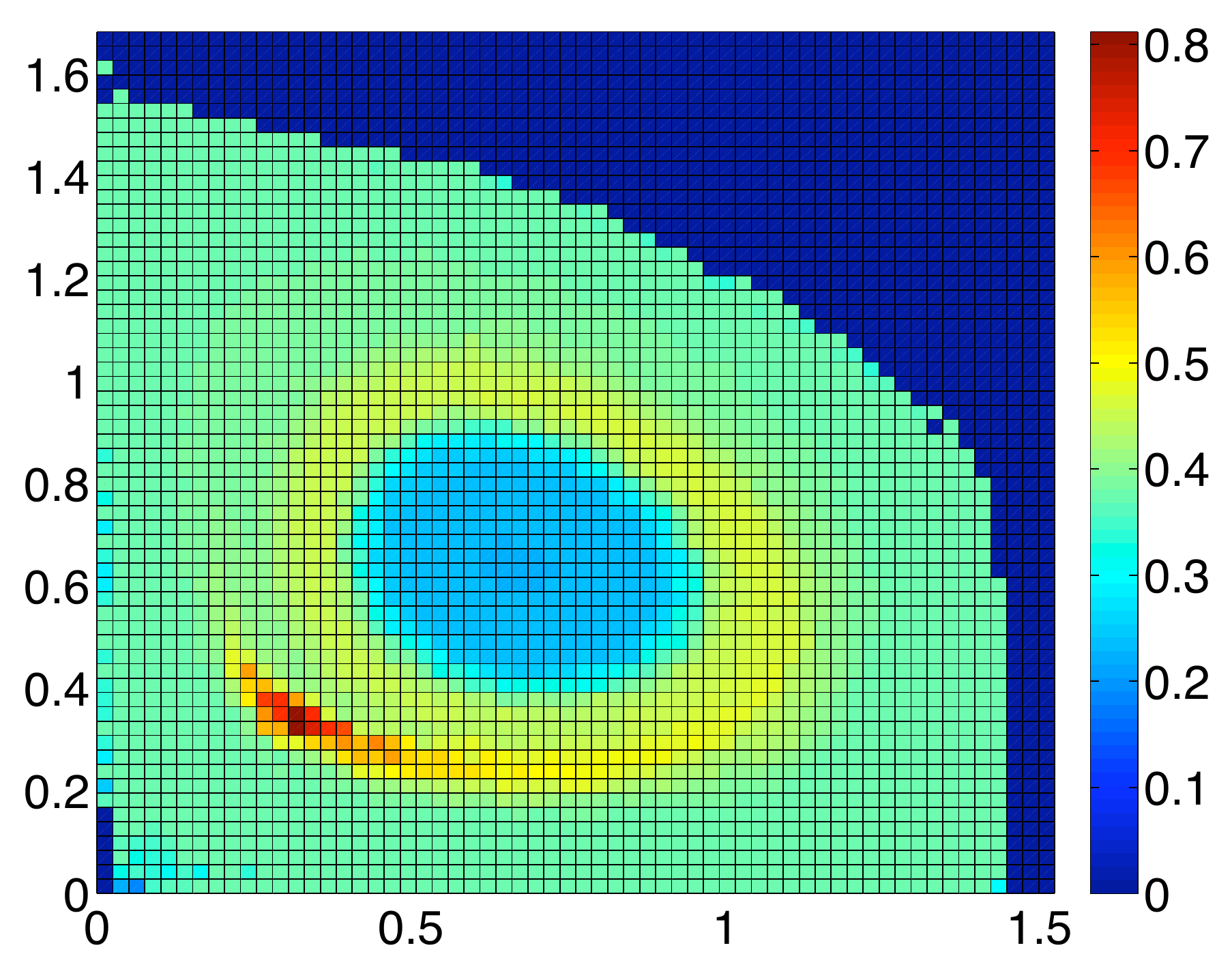}}}
\hspace{-0.1in}
\subfigure[]{\scalebox{.35}{\includegraphics{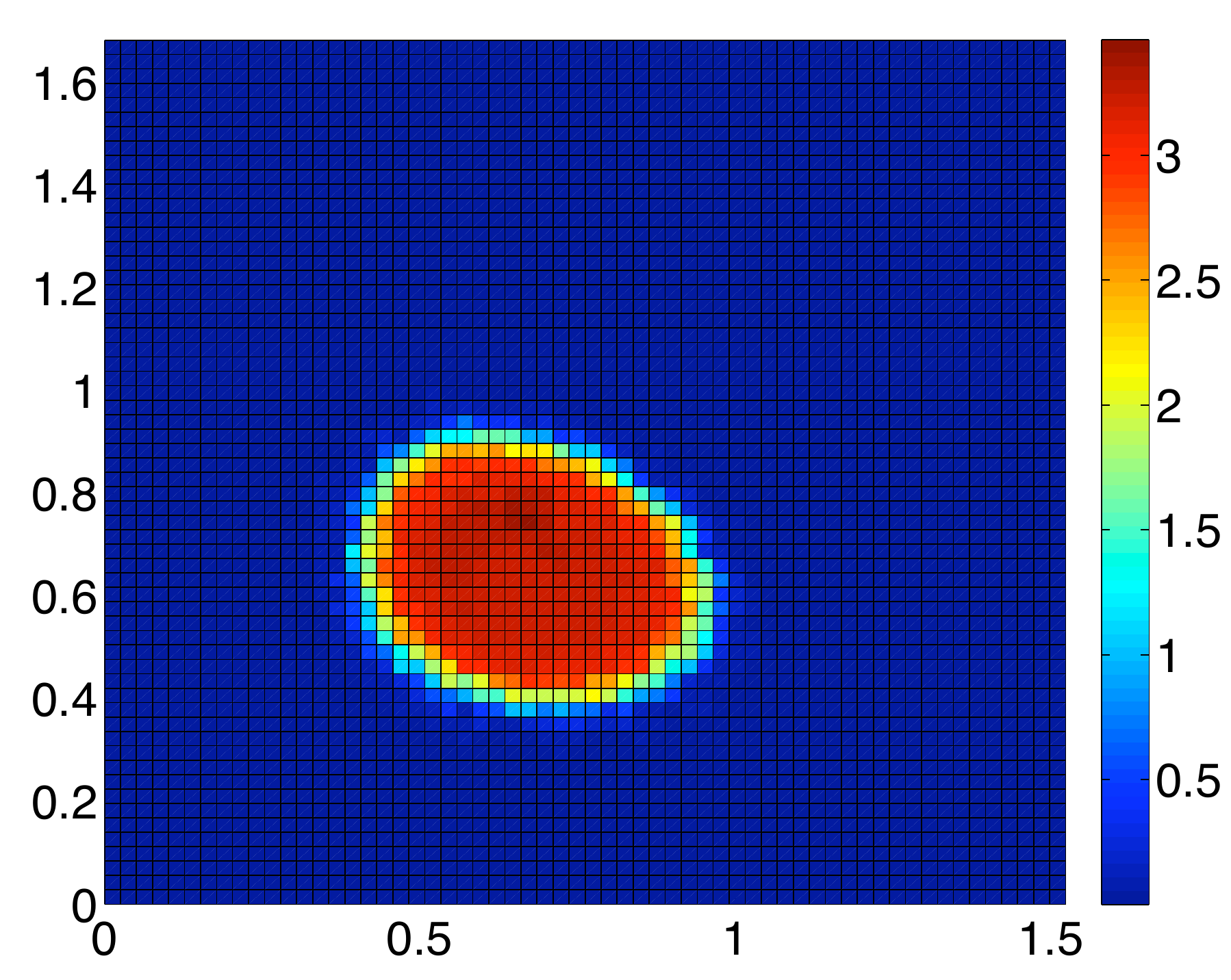}}} }
\vspace{-0.1in}
\caption{Observability gramian after 1000 iterations with sensor location at a) $B_2$; b) $B_3$. }
\label{observe23}
\end{center}
\end{figure}
%\begin{figure}[H]\centering \subfigure[]
%{\includegraphics[width=2.0 in]{3a}} \subfigure[]
%{\includegraphics[width=2.0 in]{3g}}
% \caption{  \label{room3}  (a)  A particular Sensor Position ($B_2$) (b) Gramian after time 100}
%\end{figure}
%
%This vector field captures the dynamics of room temperature fluctuation. The vector field was obtained as a steady state solution of Navier-Stokes equation. The vector field is shown in Figure \ref {vec}.
In Fig. \ref{control12}b and Fig. \ref{control3}a, we show the plots for the support of the controllability gramian after 10000 time steps corresponding to two different locations of actuator sets $B_2$ and $B_3$ respectively.  The support of the controllability gramian corresponding to $B_2$ and $B_3$  locations of actuator sets is approximately the same and equals $1.6$. However the $2$-norm of the gramian corresponding to actuator location on set $B_2$ is equal to $38$, while for $B_3$ it is equal to $35$. Comparing figures \ref{control12}a, \ref{control12}b, and  \ref{control3}a, we see that the support of the gramian for actuator location at set $B_1$ is considerably smaller but it has considerably larger $2$-norm compared to actuators locations at $B_2$ and $B_3$.  The large value of gramian with small support in Fig. \ref{control12}a can be very effective if one desires to perform localized control action. Comparing the support and the $2$-norm of the gramian function, one can conclude that the actuator location corresponding to $B_2$ is optimal among  $B_1,B_2$, and $B_3$.

In Fig. \ref{control3}b, \ref{observe23}a, and \ref{observe23}b, we show the support and $2$-norm of observability gramian corresponding to sensor locations on sets $B_1, B_2$, and $B_3$ respectively. The large support of the observability gramian for sensor location on set $B_2$ nearly outweigh the larger value and smaller support of gramian corresponding to sensors locations on set $B_1$ and $B_3$.
\begin{rem}\label{collocation}
The actuator and sensor locations in the simulation example seem to be collocated. However, this is just a coincidence and  in general this may not be the case. However it will be interesting question for future investigation. In particular, the combined problem of sensor and actuator placement will be the topic of our future investigation.

% actuator and sensor locations need not be collocated. The reason is because, following the lines of Theorem (\ref{claimcontrollabilitygramian}), one can show that the observability gramian
%$\mathcal{O}_A^{\tau} = \int_0^{\tau} (\mathbb{U}_s \chi_A(x))ds$ for an
%observer situated in the set $A$ is invertible in the following region $\mathcal{R}_A^{\tau} = \cup_{t=0}^{\tau} \phi_{-t}(B)$. This set $\mathcal{R}_A^{\tau}$
%is the set of points that trajectories flow through before ending up in the set $A$ in time $\tau$. Depending on the
%flow $f(x)$, there could be situations where $\mathcal{R}_A^{\tau}$
%is extremely small compared to the support of the controllability gramian $\mathcal{R}_B^{\tau}$, or vice versa. In the simulation, we just stumbled upon a
%set which was optimal for both actuation and sensing.
\end{rem}
%For the set $B_1$ the Lebesgue measure of the support of $\rho_{B_1}^{\tau}$ is $0.4716$ and it has an $L^2$ norm of $32.78$ at time $\tau =100$. The Actuator position and gramian plots are shown in Figure \ref {room1} (a) and (b) respectively.

%Next, we consider the sensor position $B_2$. The Lebesgue measure of the support of $\rho_{B_2}^{\tau}$ is $1.6264$ and it has an $L^2$ norm of $3.7126$ at time $\tau =100$. The Actuator position and gramian plots are shown  in Figure \ref {room2} (a) and (b) respectively.

%Finally, the sensor location is chosen as $B_3$. The Lebesgue measure of the support of $\rho_{B_3}^{\tau}$ is $1.6264$ and it has an $L^2$ norm of $6.0796$ at time $\tau =100$. The Actuator position and gramian plots are shown  in Figure \ref{room3} (a) and (b) respectively.

%Figures \ref {room1}, {room2}, and {room3} shows three different locations of actuator set (namely $B1, B2, B3$) and their corresponding gramian $\rho_{B_i}^{\tau}$ after time $\tau =100 $. The first criterion of selection for the actuator set would be maximum support for the set $\rho_{B}^{\tau}$. For $B_1$, the Lebesgue measure for the support of the finite time gramian is smallest among the three. According to criterion (1), $B_1$ would be the worst possible choice of actuator. Among $B_2$ and $B_3$ no one is better than the other if judged from criterion (1). The Lebesgue measure of the finite time gramian is the same in these cases. But $L^2$ norm for the gramian is higher for $B_3$. Thus $B_3$ would be the best possible choice among the three if we use criterion (2).
\section{Conclusion}\label{conclusion}
In this paper, controllability and observability gramian based test criteria are used to decide the suitability of given actuator and sensor locations. As compared to purely computational based methods currently existing in the literature, our proposed approach provides a systematic and insightful method for deciding the location of actuators and sensors in building systems. In particular,
stability properties of the advection vector field are shown to play an important role in deciding the location of actuators and sensors. In our future research work, the explicit formula for the gramians will be exploited to provide a systematic algorithm for determining the optimal location of sensors and actuators. Furthermore some of the assumptions made in the derivation of control equations will be removed by incorporating elements of  complex physics involved in building systems.
\section{Acknowledgement}
The authors would like to thank Prof. Baskar Ganapathysubramanian for providing data for the fluid flow vector field. Financial support of National Science Foundation Grant CMMI-0807666 is greatly acknowledged.

\bibliographystyle{IEEEtran}
\bibliography{ref}
% \appendix
% \section{}
% \label{}
\end{document}